\documentclass{article}
\usepackage{graphicx}
\usepackage{amsmath}

\newtheorem{theorem}{Theorem}[section]
\newtheorem{lemma}[theorem]{Lemma}

\begin{document}
\begin{center}
\begin{Large}
{\textbf Computing the spectrum of non self-adjoint\\[0pt]
Sturm-Liouville problems with\\[0pt]
parameter dependent boundary conditions}\\[0pt]
\end{Large}
\vspace*{1cm} by\\[0pt]
\vspace*{1cm} \textbf{B. Chanane}\\[0pt]\vspace*{1cm}
Mathematical Science Department,\\[0pt]
K.F.U.P.M., Dhahran 31261, Saudi Arabia\\[0pt]
E-Mail: chanane@kfupm.edu.sa\\[0pt]
\end{center}

\pagenumbering{arabic} \setcounter{page}{1}

\textbf{Abstract} --- This paper deals with the computation of the
eigenvalues of non self-adjoint Sturm-Liouville problems with parameter
dependent boundary conditions using the \textit{regularized sampling method}.

A few numerical examples among which singular ones will be presented to
illustrate the merit of the method and comparison made with the exact
eigenvalues when they are available.


\section{Introduction}

\setcounter{equation}{0} Non self-adjoint eigenvalue problems arise, as is
well known, in hydrodynamic and magnetohydrodynamic stability \cite%
{C55,C61,DR81,L88} while self-adjoint problems arise mostly in
quantum mechanics \cite{EE87}. The lack of oscillation theorems in
the non self-adjoint case makes any computation of the spectrum a
very difficult task \cite{D98}. In fact, the eigenvalues are
scattered over the complex plane and we need first to determine
the regions which contain them. A method that finds the
eigenvalues in a rectangle and in a left half plane has been
introduced in \cite{GM2001}. It is based on the argument principle
with compound matrix method using Magnus expansion. In
\cite{BLMTW03} the authors report on a method that provides bounds
for the eigenvalues of singular Sturm-Liouville problems over
$[0,\infty)$~with a complex potential. The method consists in
obtaining first a floating point approximation to the desired
eigenvalue by truncating the infinite interval then use interval
arithmetic to localize the eigenvalue. In \cite{B01}, the author
uses the sampling method introduced in \cite{BC96} to compute the
eigenvalues of non self-adjoint Sturm-Liouville problems.

For the mathematical foundation one may consult \cite{DS71,N68a,EE87}. On
the numerical side \cite{P93,HS97} summarize most of the available software
dealing with the computation of the eigenvalues of Sturm-Liouville problems.

In \cite{C05a}, this author introduced the \textit{regularized sampling
method}; a method which is based on Shannon's sampling theory but applied to
regularized functions. Hence avoiding any (multiple) integration and keeping
the number of terms in the Cardinal series manageable. It has been
demonstrated that the method is capable of delivering higher order estimates
of the eigenvalues at a very low cost. The purpose in this paper is to
extend the domain of application of this method to the problem at hand.

\section{Main results}

Consider the following non self-adjoint Sturm-Liouville problem with
non-separated parameter dependent boundary conditions,

\begin{equation}
\left\{
\begin{array}{c}
-y^{\prime \prime }+q(x)y=\mu^2 y\qquad ,\qquad x\in \lbrack 0,1 ] \\
A(y(0),y^{\prime }(0),y(1),y^{\prime }(1 ))^{T}=0%
\end{array}
\right.  \label{eq1}
\end{equation}
where the matrix $A(\mu)=\left(
\begin{array}{cccc}
a_{11}(\mu) & a_{12}(\mu) & a_{13}(\mu) & a_{14}(\mu) \\
a_{21}(\mu) & a_{22}(\mu) & a_{23}(\mu) & a_{24}(\mu)%
\end{array}
\right) $ has rank 2, and $q$ is a complex-valued function
satisfying $q\in L^{1}_{loc}(0,1 )$. We shall not make any
assumption on the analyticity of $A $~nor on the growth of its
components.

The purpose in this paper is to compute the eigenvalues of
(\ref{eq1}) with the
minimum of effort and a greater precision using the newly introduced \textit{%
regularized sampling method} \cite{C05a}, an improvement on the
method based on sampling theory introduced in \cite{BC96}. We note
here that the analyticity of $A$~and the conditions on the growth
of its components imposed in \cite{C05a} are not necessary for the
computation of the eigenvalues as shall be seen in the sequel. In
fact all what is needed is the recovery of certain entire
functions $h_{kl}$~associated with some base problems defined
below.

It is well known that the spectrum is discrete and scattered over the
complex plane which makes difficult its computation. Also, there is no
result about the distribution nor the multiplicity of the eigenvalues.

Let $y_{c}(x,\mu)$~and $y_{s}(x,\mu)$, be the solutions of the base problems
\begin{equation}
\left\{
\begin{array}{c}
-y^{\prime \prime }+q(x)y=\mu^2 y\qquad ,\qquad x\in \lbrack 0,1] \\
y(0)=1~,~y^{\prime}(0)=0%
\end{array}
\right.
\end{equation}
and
\begin{equation}
\left\{
\begin{array}{c}
-y^{\prime \prime }+q(x)y=\mu^2 y\qquad ,\qquad x\in \lbrack 0,1] \\
y(0)=0~,~y^{\prime}(0)=1%
\end{array}
\right.
\end{equation}
respectively. Then the general solution of the differential equation in (\ref%
{eq1}) and its derivative are
\begin{eqnarray*}
y(x,\mu)&=&c_{1}y_{c}(x,\mu)+c_{2}y_{s}(x,\mu) \\
y^{\prime}(x,\mu)&=&c_{1}y_{c}^{\prime}(x,\mu)+c_{2}y_{s}^{\prime}(x,\mu)
\end{eqnarray*}
The boundary condition gives after separating $c_{1} $ and $c_{2}$,
\begin{equation}
c_{1} Aw_{1}+c_{2} Aw_{2}=0  \label{eq7}
\end{equation}
where
\begin{eqnarray*}
w_{1} &=&\left( 1,0,y_{c}(1,\mu) ,y_{c}^{\prime}(1,\mu)\right) ^{T}
\label{eq8} \\
w_{2} &=&\left( 0,1,y_{s}(1,\mu) ,y_{s}^{\prime}(1,\mu) \right) ^{T}
\label{eq9}
\end{eqnarray*}
Thus, a necessary and sufficient condition for $\lambda=\mu^2$~to be an
eigenvalue is that $\mu$~satisfies the characteristic equation $B(\mu )=0$,
where $B$~is the characteristic function $B(\mu)=\det \left(
Aw_{1}|Aw_{2}\right) =\det [A\left( w_{1}|w_{2}\right)]$, that is,
\begin{eqnarray*}
B(\mu)&=&
\left(a_{11}(\mu)+a_{13}(\mu)y_{c}(1,\mu)+a_{14}(\mu)y_{c}^{\prime}(1,\mu)%
\right)
\left(a_{22}(\mu)+a_{23}(\mu)y_{s}(1,\mu)+a_{24}(\mu)y_{s}^{\prime}(1a,\mu)%
\right) \\
&&-\left(a_{21}(\mu)+a_{23}(\mu)y_{c}(1,\mu)+a_{24}(\mu)y_{c}^{\prime}(1,%
\mu)\right)
\left(a_{12}(\mu)+a_{13}(\mu)y_{s}(1,\mu)+a_{14}(\mu)y_{s}^{\prime}(1,\mu)%
\right)
\end{eqnarray*}

We shall need the following well known results,

\begin{lemma}
\cite{C97} $\sin z/z$~and $\cos z$~are entire as functions of $z$~and
satisfy the estimates
\begin{equation}
\left\vert \sin z/z\right\vert \leq \beta_{0}e^{\left\vert
Im~z\right\vert }/(1+\left\vert z\right\vert
)\mbox{~and~}\left\vert \cos z\right\vert \leq e^{\left\vert
Im~z\right\vert }  \notag
\end{equation}%
where $\beta_{0}=1.72$.
\end{lemma}

Using the above lemma one can show the following result to hold.

\begin{theorem}
\cite{C97} $y_{c}(x,\mu )$,~$y_{s}(x,\mu )$,~$y_{c}^{\prime }(x,\mu )$~and $%
y_{s}^{\prime }(x,\mu )$~are entire as functions of $\mu $~for each fixed $%
x\in (0,1]$~and satisfy the growth conditions,
\begin{equation*}
|y_{c}(x,\mu )-\cos (\mu x)|,~|y_{s}(x,\mu )-\frac{\sin (\mu
x)}{\mu }|, |y_{c}^{\prime }(x,\mu )+\mu \sin (\mu
x))|~,~|y_{s}^{\prime }(x,\mu )-\cos (\mu x))| \leq \beta
_{1}e^{x|\mbox{Im}\mu |}
\end{equation*}
for some positive constant $\beta _{1}$.
\end{theorem}

In \cite{C97} and \cite{C2001} we have obtained much higher
estimates of the eigenvalues than those presented in Theorem 2.2
above, at the expense of subtracting terms involving multiple
integrals. Here and as in \cite{C05a}, we shall stick with the
estimates given in Theorem 2.2, avoiding any (multiple)
integration. We shall show by the same token that we can get a
higher order estimate of the eigenvalues of the problem at hand at
a very low cost. In fact we do not have even to keep on increasing
the number of sampling points.

Let $PW_{\sigma }$~denote the Paley-Wiener space \cite{Z93}
\begin{equation*}
PW_{\sigma }=\{f\mbox{entire, }|f(\mu )|\leq \beta\mbox{e}^{\sigma
|\mbox{Im}\mu |}\mbox{, }\int_{\mbox{R}}|f(\mu )|^{2}d\mu <\infty
\}
\end{equation*}

Let $h_{kl}$~be defined by
\begin{equation}
\left\{
\begin{array}{ccc}
h_{11}(\mu) & = & \left( \frac{\sin \theta \mu }{\theta \mu }\right)
^{m}\left(y_{c}(1,\mu)-\cos \mu\right) \\
h_{12}(\mu) & = & \left( \frac{\sin \theta \mu }{\theta \mu }\right)
^{m}\left(y_{s}(1,\mu)-\frac{\sin \mu }{\mu}\right) \\
h_{21}(\mu) & = & \left( \frac{\sin \theta \mu }{\theta \mu }\right)
^{m}\left(y_{c}^{\prime}(1,\mu)+\mu\sin \mu \right) \\
h_{22}(\mu) & = & \left( \frac{\sin \theta \mu }{\theta \mu }\right)
^{m}\left(y_{s}^{\prime}(1,\mu)-\cos \mu\right)%
\end{array}
\right.  \notag
\end{equation}

Then we rewrite $y_{c}(1,\mu),~y_{c}^{\prime}(1,\mu),~y_{s}(1,\mu)$~and $%
y_{s}^{\prime}(1,\mu)$~as
\begin{equation}
\left\{
\begin{array}{ccc}
y_{c}(1,\mu) & = & h_{11}(\mu)\left( \frac{\sin \theta \mu }{\theta \mu }%
\right) ^{-m}+\cos \mu \\
y_{s}(1,\mu) & = & h_{12}(\mu)\left( \frac{\sin \theta \mu }{\theta \mu }%
\right) ^{-m}+\frac{\sin \mu }{\mu} \\
y_{c}^{\prime}(1,\mu) & = & h_{21}(\mu)\left( \frac{\sin \theta \mu }{\theta
\mu }\right) ^{-m}-\mu\sin \mu \\
y_{s}^{\prime}(1,\mu) & = & h_{22}(\mu)\left( \frac{\sin \theta \mu }{\theta
\mu }\right) ^{-m}+\cos \mu%
\end{array}
\right.  \notag
\end{equation}

\begin{theorem}
Let $\vartheta$~be a positive constant and $m$~ be a positive integer ($%
m\geq 2$). The functions $h_{kl},~(k,l=1,2)$~ belong to the Paley space $%
PW_{\sigma }$~with $\sigma=1+m\theta$ and satisfy the estimates
\begin{equation}
\left| h_{kl}(\mu)\right| \leq \frac{\beta_{2}}{(1+\theta |\mu |)^{m}}%
e^{\sigma\left| \mbox{Im}\mu \right| }  \notag
\end{equation}
$k,l=1,2$~for some positive constant $\beta_{2}$.
\end{theorem}

\textbf{Proof:} That $h_{kl}$~are entire and satisfy the given
estimates is a direct consequence of Theorem 2.2 and the fact that
$\frac{\sin \theta \mu }{\theta \mu }$ is an entire function of
$\mu $ and satisfy the estimate in Lemma 2.1$\diamondsuit$

Since the $h_{kl}(\mu)$ belong to the Paley-Wiener space $PW_{\sigma }$ ~for
each $k,l=1,2$, they can be recovered from their values at the points\ $\mu
_{j}=j\frac{\pi }{\sigma }$, $j\in Z$, using the following celebrated
theorem,

\begin{theorem}[Whitaker-Shannon-Kotel'nikov]
\cite{Z93} Let $h\in PW_{\sigma }$, then
\begin{equation}
h(\mu )=\sum_{j=-\infty }^{\infty }h(\mu _{j})\frac{\sin \sigma (\mu -\mu
_{j})}{\sigma (\mu -\mu _{j})}  \notag
\end{equation}%
$\mu _{j}=j\frac{\pi }{\sigma }$. The series converges absolutely and
uniformly on compact subsets of $C$ and in $L_{d\mu }^{2}(R)$.
\end{theorem}

For all practical purposes, we consider finite summations, therefore we need
to approximate $h_{kl}$~by a truncated series $h_{kl}^{[N]}$. The following
lemma gives an estimate for the truncation error.

\begin{lemma}[Truncation error]
Let $h_{kl}^{[N]}(\mu )=\sum_{j=-N}^{N}h_{kl}(\mu _{j})\frac{\sin \sigma
(\mu -\mu _{j})}{\sigma (\mu -\mu _{j})}$ denote the truncation of $%
h_{kl}(\mu )$. Then, for $|\mu |<N\pi /\sigma $,

\begin{equation*}
\left| h_{kl}(\mu )-h_{kl}^{[N]}(\mu )\right| \quad \leq
\frac{|\sin \mu |\beta_{3}}{\pi (\pi /\sigma
)^{m-1}\sqrt{1-4^{-m+1}} }\left[ \frac{1}{\sqrt{(N\pi /\sigma
)-\mu }}+\frac{1}{\sqrt{(N\pi /\sigma )+\mu }}\right]
\frac{1}{(N+1)^{m-1}},
\end{equation*}
where $\beta_{3}=||\mu ^{m-1}h_{kl}(\mu )||_{2}$.
\end{lemma}

\textbf{Proof:} Since $\mu ^{m-1}h_{kl}(\mu )\in L^{2}(-\infty
,\infty )$, Jagerman's result (see \cite{Z93}, Theorem 3.21, p.90)
is applicable and yields the given estimate for the
$h_{kl},~k,l=1,2$$\diamondsuit$

An approximation $B_{N}$~to the characteristic function $B$~is provided by
replacing the $h_{kl}$~by its approximation $h_{kl}^{[N]}$, and we obtain at
once,

\begin{lemma}
The approximate characteristic function $B_{N}$~satisfies the estimate,

\begin{equation*}
\left| B(\mu )-B_{N}(\mu )\right| \quad \leq \left|\frac{\sin
\theta \mu }{\theta \mu }\right|^{-m}\frac{ |\sin \mu
|\beta_{4}}{\pi (\pi /\sigma )^{m-1}\sqrt{1-4^{-m+1}}}\left[
\frac{1}{\sqrt{(N\pi /\sigma )-\mu }}+\frac{1}{\sqrt{(N\pi /\sigma
)+\mu }}\right] \frac{1}{(N+1)^{m-1}},
\end{equation*}
for some positive constant $\beta_{4}$.
\end{lemma}

We claim the following,

\begin{theorem}
Let $\overline{\mu}^{2}$~be an exact eigenvalue of $B$~ of multiplicity $n$
and denote by $\mu_{N}^{2}$~the corresponding approximation of a square of a
zero of $B_{N}$. Then, for $|\mu_{N}|<N\pi/\sigma$, we have,
\begin{eqnarray*}
|\mu_{N}-\overline{\mu}|&\leq&\left(\frac{m!}{\inf|B^{(m)}(\widetilde{\mu})|}%
\left|\frac{\sin \theta \mu_{N} }{\theta \mu_{N } }\right|^{-m}\frac{|\sin
\mu_{N }|\beta_{4}}{\pi (\pi /\sigma )^{m-1}\sqrt{1-4^{-m+1}}}\right)^{1/n} \\
&&\times\left[ \frac{1}{\sqrt{(N\pi /\sigma )-\mu_{N } }}+\frac{1}{\sqrt{%
(N\pi /\sigma )+\mu_{N } }}\right]^{1/n} \frac{1}{(N+1)^{(m-1)/n}}
\end{eqnarray*}
where the $\inf$~is taken over a ball centered at $\mu_{N}$~with radius $%
|\mu_{N}-\overline{\mu}|$~and not containing a multiple of $\pi/\theta$.
\end{theorem}

\textbf{Proof} Since $\overline{\mu}$~is a zero of $B$~with multiplicity $n$%
, then
\begin{equation}
B(\overline{\mu })-B(\mu _{N})=\frac{(\overline{\mu}-\mu_{N})^n}{n!}B^{(n)}(%
\widetilde{\mu})  \notag
\end{equation}
for some $\widetilde{\mu}$. Thus,
\begin{eqnarray*}
|\overline{\mu}-\mu_{N}|^n&=& \frac{m!|B(\overline{\mu })-B(\mu _{N})|}{
|B^{(m)}(\widetilde{\mu})|} \\
&\leq&\frac{m!}{\inf |B^{(m)}(\widetilde{\mu})|}\left|\frac{\sin \theta \mu_{N} }{%
\theta \mu_{N } }\right|^{-m}\frac{|\sin \mu_{N }|\beta_{4}}{\pi
(\pi /\sigma
)^{m-1}\sqrt{1-4^{-m+1}}} \\
&&\times\left[ \frac{1}{\sqrt{(N\pi /\sigma )-\mu_{N } }}+\frac{1}{\sqrt{%
(N\pi /\sigma )+\mu_{N } }}\right] \frac{1}{(N+1)^{m-1}}
\end{eqnarray*}
where the $\inf$~is taken over a ball centered at $\mu_{N}$~with radius $%
|\mu_{N}-\overline{\mu}|$~and not containing a multiple of
$\pi/\theta$. Thus, the result$\diamondsuit$

\section{Numerical examples}

\setcounter{equation}{0} In this section, we shall present a few examples to
illustrate our method. We have taken $\theta =1/(N-m)$ in order to avoid the
first singularity of $\left( \frac{\sin \theta \mu _{N}}{\theta \mu _{N}}%
\right) ^{-1}$. The sampling values were obtained using the Fehlberg 4-5
order Runge-Kutta method. The first two problems are taken from \cite%
{BLMTW03} in which the authors use interval arithmetic to localize
the eigenvalues of singular Sturm-Liouville problems with complex
potentials. The third problem, taken from \cite{B01}, shows that
the \textit{regularized sampling method} provides much better
results than the sampling method without regularization. The last
example demonstrates that our method can estimate the eigenvalues
with a great precision even in situation where other methods might
introduce spurious eigenvalues and/or miss some of them. We shall
mention however that we shall not make use of the error estimate
given above for the time being. The method consists first in the
recovery of the entire functions $h_{kl}$ with great precision,
then use the boundary conditions to determine the characteristic
function. The zeros of this characteristic function are the square
roots of the sought eigenvalues. We shall denote
$\iota=\sqrt{-1}$.

\textbf{Example 3.1}(Taken from \cite{BLMTW03}) Consider the
singular Sturm-Liouville problem
\begin{equation}
\left\{
\begin{array}{c}
-y^{\prime \prime }(x)+10~\iota\sin x e^{-x}y(x)=\lambda
y(x)~,~0\leq x<\infty ,
\\
y(0)=0%
\end{array}
\right.  \notag
\end{equation}

We shall use interval truncation and compute the eigenvalues of
\begin{equation}
\left\{
\begin{array}{c}
-y_{\gamma}^{\prime \prime }(x)+10~\iota\sin x
e^{-x}y_{\gamma}(x)=\mu^2
y_{\gamma}(x)~,~0\leq x\leq \gamma , \\
y_{\gamma}(0)=0~,~y_{\gamma}^{\prime}(\gamma)=~\iota\mu y_{\gamma}(\gamma)%
\end{array}
\right.  \notag
\end{equation}
and as in \cite{BLMTW03} we shall take $\gamma=10$. The second boundary
condition has been obtained by considering the Jost solution $y=e^{\iota \mu x}$%
~and its derivative $y^{\prime}=\iota\mu e^{\iota\mu x}$, thus, $y^{\prime}(%
\gamma)=\iota\mu~y(\gamma)$. In \cite{BLMTW03} the authors
obtained an eigenvalue lying in $1.6043912_{58}^{64} +
1.7978849_{67}^{81}~\iota$~where the notation $2.1_{6}^{4}$~stands
for the interval $[2.14,2.16]$. Taking $N=40$, and for different
values of $m$~we obtained the results summarized in Table 1,

\begin{center}
\begin{equation*}
\begin{tabular}{|c|c|}
\hline $m$ & Approximate Eigenvalue \\
\hline 5 & 1.604391348283 + 1.797884747658~$\iota$ \\
\hline 10 & 1.604391251270 + 1.797884973775~$\iota$\\
 \hline 15 & 1.604391251323 + 1.797884973746~$\iota$ \\
 \hline
\end{tabular}\label{tab2a}
\end{equation*}
Table 1: Approximation of an eigenvalue for different values of
$m$ in Example 3.1
\end{center}

\textbf{Example 3.2}(Taken from \cite{BLMTW03}) Consider the
singular problem
\begin{equation}
\left\{
\begin{array}{c}
-y^{\prime \prime }(x)+10~\iota~e^{-x}y(x)=\lambda y(x)~,~0\leq x<\infty , \\
y(0)=0%
\end{array}
\right.  \notag
\end{equation}

We shall use interval truncation and compute the eigenvalues of
\begin{equation}
\left\{
\begin{array}{c}
-y_{\gamma}^{\prime \prime }(x)+10~\iota e^{-x}y_{\gamma}(x)=\mu^2
y_{\gamma}(x)~,~0\leq x\leq \gamma , \\
y_{\gamma}(0)=0~,~y_{\gamma}^{\prime}(\gamma)=~\iota\mu y_{\gamma}(\gamma)%
\end{array}
\right.  \notag
\end{equation}
and as in \cite{BLMTW03} we shall take $\gamma=10$. In
\cite{BLMTW03} the authors obtained an eigenvalue lying in
$2.8122672_{89}^{92} + 2.17223818_{78}^{99}~\iota$. Taking $N=40$,
and for different values of $m$~we obtained the results summarized
in Table 2,

\begin{center}
\begin{equation*}
\begin{tabular}{|c|c|}
\hline $m$ & Approximate Eigenvalue \\ \hline 5 &
2.812264032443898167911 + 2.17223723666731852353~$\iota$ \\ \hline
10 & 2.8122672894628469454261 + 2.172238191264223861~$\iota$ \\
\hline 15 & 2.812267288417814133626 + 2.172238191179093864~$\iota$ \\
\hline
\end{tabular}%
\label{tab2b}%
\end{equation*}
Table 2: Approximation of an eigenvalue for different values of
$m$ in Example 3.2

\end{center}

\textbf{Example 3.3}(Taken from \cite{B01}) Consider the non
self-adjoint problem
\begin{equation}
\left\{
\begin{array}{c}
-y^{\prime \prime }(x)+(3-2~\iota)y(x)=\lambda y(x)~,~0\leq x\leq \pi , \\
y(0)=y(\pi)=0%
\end{array}
\right.  \notag
\end{equation}

The exact eigenvalues of the original problem are
$\lambda_{k}=k^2+3-2~\iota,~k=1,2,...$. Taking $N=40$~and $m=10$,
we obtained the results summarized in Table 3,

\begin{center}
\begin{equation*}
\begin{tabular}{|c|c|c|c|}
\hline Index & Exact Eigenvalue & Approximate Eigenvalue &
Absolute Error \\ \hline 1 & 4-2~$\iota$ &
3.9999999999999289-1.99999999999998304519~$\iota$ & 7.30$\times
10^{-14}$
\\ \hline
2 & 7-2~$\iota$ &
6.9999999999998187-1.99999999999984438039~$\iota$ & 2.38$\times
10^{-13}$
\\ \hline
3 & 12-2~$\iota$ &
11.999999999999561-1.9999999999997718164~$\iota$ & 4.93$\times
10^{-13}$
\\ \hline
4 & 19-2~$\iota$ &
18.999999999999172-2.0000000000000542265~$\iota$ & 8.29$\times
10^{-13}$
\\ \hline
5 & 28-2~$\iota$ &
27.999999999999391-2.0000000000011521175~$\iota$ & 1.30$\times
10^{-12}$
\\ \hline
6 & 39-2~$\iota$ &
39.000000000001586-2.0000000000029645542~$\iota$ & 3.36$\times
10^{-12}$
\\ \hline
7 & 52-2~$\iota$ &
52.000000000005729-2.0000000000033954538~$\iota$ & 6.66$\times
10^{-12}$
\\ \hline
8 & 67-2~$\iota$ &
67.000000000006628-1.9999999999977990747~$\iota$ & 6.98$\times
10^{-12}$
\\ \hline
9 & 84-2~$\iota$ &
83.999999999993226-1.9999999999829522498~$\iota$ & 1.83$\times
10^{-11}$
\\ \hline
10 & 103-2~$\iota$ &
102.999999999961-1.9999999999690674138~$\iota$ & 4.95$\times
10^{-11}$
\\ \hline
11 & 124-2~$\iota$ &
123.999999999944-1.9999999999941154414~$\iota$ & 5.55$\times
10^{-11}$
\\ \hline
12 & 147-2~$\iota$ &
147.000000000038-2.0000000001107296594~$\iota$ & 1.17$\times
10^{-10}$
\\ \hline
13 & 172-2~$\iota$ &
172.000000000323-2.0000000002862574821~$\iota$ & 4.32$\times
10^{-10}$
\\ \hline
14 & 199-2~$\iota$ &
199.000000000556-2.0000000001798594170~$\iota$ & 5.85$\times
10^{-10}$
\\ \hline
15 & 228-2~$\iota$ &
227.999999999678-1.9999999989989010834~$\iota$ & 1.05$\times
10^{-9}$
\\ \hline
16 & 259-2~$\iota$ &
258.999999996079-1.9999999963552374114~$\iota$ & 5.35$\times
10^{-9}$
\\ \hline
17 & 292-2~$\iota$ &
291.999999991669-1.9999999961982129591~$\iota$ & 9.15$\times
10^{-9}$
\\ \hline
18 & 327-2~$\iota$ &
327.000000004474-2.0000000150098733387~$\iota$ & 1.56$\times
10^{-8}$
\\ \hline
19 & 364-2~$\iota$ &
364.000000082618-2.0000000798746168342~$\iota$ & 1.14$\times
10^{-7}$
\\ \hline
20 & 403-2~$\iota$ &
403.000000232131-2.0000001291033812382~$\iota$ & 2.65$\times
10^{-7}$
\\ \hline
\end{tabular}%
\label{tab4}%
\end{equation*}
Table 3: Exact and Approximate eigenvalues in Example 3.3

\end{center}

\textbf{Example 3.4} Consider now the following non self-adjoint
Sturm-Liouville problem with complex potential and parameter dependent
boundary condition,
\begin{equation}
\left\{
\begin{array}{c}
-y^{\prime \prime }(x)+e^{2~\iota x}y(x)=\mu^2 y(x)~,~0\leq x\leq 1 , \\
y(0)+\mu y(1)=0 \\
y^{\prime}(0)=0%
\end{array}
\right.  \notag
\end{equation}
Here again we are in a position to derive the exact characteristic function
which in fact can be expressed in terms of Bessel functions. Indeed, let $%
\lambda =\mu ^{2}$ and consider the change of variables $t=
\mathbf{e}^{\iota x}$. The differential equation becomes the
Bessel equation of order $\mu $ given by
\begin{equation}
t^{2}\frac{d^{2}z}{dt^{2}}+t\frac{dz}{dt}+(t^{2}-\mu ^{2})z=0  \notag
\end{equation}
whose solution is
\begin{equation}
z(t)=c_{1}\mathcal{J}_{\mu }(t)+c_{2}\mathcal{J}_{-\mu }(t),  \notag
\end{equation}
where $\mathcal{J}_{\mu }$ and $\mathcal{J}_{-\mu }$ are the Bessel
functions of the first kind of order $\mu $.

Returning to the original variables, we obtain
\begin{equation}
y(x)=c_{1}\mathcal{J}_{\mu }(\mathbf{e}^{\iota x})+c_{2}\mathcal{J} _{-\mu }(%
\mathbf{e}^{\iota x}).  \notag
\end{equation}
Taking into account the boundary conditions, we obtain the homogeneous
system in $c_{1}$and $c_{2}$
\begin{equation}
\left\{
\begin{array}{c}
c_{1}\mathcal{J}_{\mu }(1)+c_{2}\mathcal{J}_{-\mu }(1)+ \mu \left( c_{1}%
\mathcal{J}_{\mu }(\mathbf{e}^{\iota})+c_{2}\mathcal{J}_{-\mu }(\mathbf{e}%
^{\iota})\right)=0 \\
c_{1}\mathcal{J}^{\prime}_{\mu }(1)+c_{2}\mathcal{J}^{\prime}_{-\mu }(1)=0%
\end{array}
\right.  \notag
\end{equation}
In order to have a nontrivial solution, a necessary and sufficient condition
is to have $B_{exact}(\mu )=0$ where

\begin{eqnarray*}
B_{exact}(\mu ) = det~\left(\!\!
\begin{array}{cc}
\mathcal{J}_{\mu }(1)+ \mu\mathcal{J}_{\mu }(\mathbf{e}^{\iota}) & \mathcal{J%
}_{-\mu }(1)+ \mu\mathcal{J}_{-\mu }(\mathbf{e}^{\iota}) \\
\mathcal{J}^{\prime}_{\mu }(1) & \mathcal{J}^{\prime}_{-\mu }(1)%
\end{array}%
\!\! \right)
\end{eqnarray*}
\unskip is the characteristic function. Now, using the well known result
\begin{equation}
\frac{d}{dx}\mathcal{J}_{\mu }(x) =\left(\mathcal{J}_{-\mu -1}(x)-\mathcal{J}%
_{\mu +1}(x)\right)/2,  \notag
\end{equation}
we obtain
\begin{eqnarray*}
B_{exact}(\mu ) = det~\left(\!\!
\begin{array}{cc}
\mathcal{J}_{\mu }(1)+ \mu\mathcal{J}_{\mu }(\mathbf{e}^{\iota}) & \mathcal{J%
}_{-\mu }(1)+ \mu\mathcal{J}_{-\mu }(\mathbf{e}^{\iota}) \\
\left(\mathcal{J}_{-\mu -1}(1)-\mathcal{J}_{\mu +1}(1)\right)/2 & \left(%
\mathcal{J}_{\mu -1}(1)-\mathcal{J}_{-\mu +1}(1)\right)/2%
\end{array}%
\!\! \right)
\end{eqnarray*}

Taking $N=40$, and $m=10$~we obtained the results summarized in
Table 4.
\begin{scriptsize}
\begin{center}
\begin{equation*}
\begin{tabular}{|c|c|c|c|c|}
\hline Index & Exact Eigenvalue & Approximate Eigenvalue &
Absolute Error & Relative Error \\ \hline 1 &
4.9685430929323576+0.3906545895360696~$\iota$ &
4.9685430929323625+0.3906545895360721~$\iota$ & 5.549$\times 10^{-15}$ & 1.113$%
\times 10^{-15}$ \\
2 & 20.60271034889337+0.75023252353154~$\iota$ &
20.60271034889340+0.75023252353155~$\iota$ & 3.393$\times
10^{-14}$ & 1.645$\times
10^{-15} $ \\
3 & 64.14038244804547+0.68422837531133~$\iota$ &
64.14038244804526+0.68422837531099~$\iota$ & 3.977$\times
10^{-13}$ & 6.201$\times
10^{-15} $ \\
4 & 119.34792168887388+0.71497240479401~$\iota$ &
119.34792168887345+0.71497240479334~$\iota$ & 8.004$\times 10^{-13}$ & 6.706$%
\times 10^{-15} $ \\
5 & 202.31443747778734+0.70057212586525~$\iota$ &
202.31443747778739+0.70057212586545~$\iota$ & 2.064$\times 10^{-13}$ & 1.020$%
\times 10^{-15} $ \\
6 & 419.44558800598641+0.70446189520144~$\iota$ &
419.44558800598892+0.70446189520528~$\iota$ & 4.582$\times 10^{-12}$ & 1.092$%
\times 10^{-14}$ \\
7 & 553.61789373762934+0.70954623577257~$\iota$ &
553.61789373762976+0.70954623577282~$\iota$ & 4.969$\times 10^{-13}$ & 8.977$%
\times 10^{-16} $ \\
8 & 715.53365857906959+0.70595783818772~$\iota$ &
715.53365857906140+0.70595783817453~$\iota$ & 1.553$\times 10^{-11}$ & 2.170$%
\times 10^{-14} $ \\
9 & 889.18520034251622+0.70898948206981~$\iota$ &
889.18520034250143+0.70898948204681~$\iota$ & 2.734$\times 10^{-11}$ & 3.075$%
\times 10^{-14} $ \\
10 & 1090.57859485902126+0.70668585309098~$\iota$ &
1090.57859485902214+0.70668585309385~$\iota$ & 3.00$\times 10^{-12}$ & 2.751$%
\times 10^{-15} $ \\
11 & 1303.70898166607058+0.70869788000992~$\iota$ &
1303.70898166611992+0.70869788008925~$\iota$ & 9.341$\times 10^{-11}$ & 7.165$%
\times 10^{-14} $ \\
12 & 1544.58037965386611+0.70709389168016~$\iota$ &
1544.58037965396658+0.70709389183628~$\iota$ & 1.856$\times 10^{-10}$ & 1.202$%
\times 10^{-13} $ \\
13 & 1797.18943505543540+0.70852627026801~$\iota$ &
1797.18943505544546+0.70852627027458~$\iota$ & 1.201$\times 10^{-11}$ & 6.687$%
\times 10^{-15} $ \\
14 & 2077.53900632820814+0.70734525957323~$\iota$ &
2077.53900632774381+0.70734525883073~$\iota$ & 8.757$\times 10^{-10}$ & 4.215$%
\times 10^{-13} $ \\
15 & 2369.6266391592291618+0.70841680475450~$\iota$ &
2369.62663915816209+0.70841680308871~$\iota$ & 1.978$\times 10^{-9}$ & 8.348$%
\times 10^{-13}$ \\
16 & 2689.45447190894724+0.70751097714777~$\iota$ &
2689.45447190899851+0.70751097732396~$\iota$ & 1.834$\times 10^{-10}$ & 6.822$%
\times 10^{-14} $ \\
17 & 3021.02063035583927+0.70834272705249~$\iota$ &
3021.02063036245257+0.70834273761505~$\iota$ & 1.246$\times 10^{-8}$ & 4.125$%
\times 10^{-12} $ \\
18 & 3380.32677490847313+0.70762595333135~$\iota$ &
3380.32677492777072+0.70762598345876~$\iota$ & 3.577$\times 10^{-8}$ & 1.058$%
\times 10^{-11} $ \\
19 & 3751.37142735725052+0.70829027475581~$\iota$ &
3751.37142735927201+0.70829027606949~$\iota$ & 2.410$\times 10^{-9}$ & 6.426$%
\times 10^{-13}$ \\
20 & 4150.15591451714336+0.70770896968510~$\iota$ &
4150.15591430123958+0.70770862603678~$\iota$ & 4.058$\times 10^{-7}$ & 9.778$%
\times 10^{-11} $ \\
21 & 4560.67904058883973+0.70825178027457~$\iota$ &
4560.67903973138474+0.70825043937887~$\iota$ & 1.591$\times 10^{-6}$ & 3.489$%
\times 10^{-10} $ \\
22 & 4998.94189026423779+0.70777085938814~$\iota$ &
4998.94189032279592+0.70777105963898~$\iota$ & 2.086$\times 10^{-7}$ & 4.173$%
\times 10^{-11} $ \\
23 & 5448.94347623327640+0.70822269647649~$\iota$ &
5448.94349859286012+0.70825825178573~$\iota$ & 0.00004200 &
7.708$\times 10^{-9} $
\\
24 & 5926.68470186115217+0.70781822661900~$\iota$ &
5926.68487236793343+0.70808488449668~$\iota$ & 0.0003165 &
5.340$\times 10^{-8}$
\\
25 & 6416.16473814590617+0.70820018792052~$\iota$ &
6416.16478405538947+0.70823000684001~$\iota$ & 0.00005474 &
8.532$\times 10^{-9} $
\\ \hline
\end{tabular}%
\label{tab5}
\end{equation*}
Table 4: Exact and Approximate eigenvalues in Example 3.4

\end{center}
\end{scriptsize}
\section{Conclusion}

In this paper, we have used the \textit{regularized sampling
method} introduced recently \cite{C05a} to compute the eigenvalues
of non self-adjoint Sturm-Liouville problems with nonseparable
parameter dependent boundary conditions. We recall that this
method constitutes an improvement upon the method based on
Shannon's sampling theory introduced in \cite{BC96} since it uses
a regularization avoiding any multiple integration. The method
allows us to get higher order estimates of the eigenvalues at a
very low cost. We have presented a few examples, including
singular ones, to illustrate the method and compared the computed
eigenvalues with the exact ones when they are available.

\section*{Acknowledgments}

The author wishes to thank King Fahd University of Petroleum and
Minerals for its constant support, and by making this possible
through grant MS/SPECTRAL/269. The thorough refereing process and
the pertinent suggestions received are greatly appreciated.

\end{document}